\documentclass[journal]{IEEEtran}

\ifCLASSINFOpdf
  \usepackage[pdftex]{graphicx}
\else
\fi

\usepackage{amsmath,bm}

\usepackage{float}
\floatstyle{ruled}
\newfloat{model}{thp}{lop}
\floatname{model}{Model}

\usepackage{xcolor}

\begin{document}

\title{The Impacts of Convex Piecewise Linear Cost Formulations on AC Optimal Power Flow}

\author{
    Carleton~Coffrin, 
    Bernard~Knueven,
    Jesse~Holzer,
    and~Marc~Vuffray
\thanks{C. Coffrin and M. Vuffray are staff scientists with Los Alamos National Laboratory, Los Alamos,
NM, 87545 USA. e-mail: cjc@lanl.gov}
\thanks{B. Knueven is a staff scientist with the National Renewable Energy Laboratory.}
\thanks{J. Holzer is a staff scientist with the Pacific Northwest National Laboratory.}
}

\markboth{}{}

\maketitle

\begin{abstract}
Despite strong connections through shared application areas, research efforts on power market optimization (e.g., unit commitment) and power network optimization (e.g., optimal power flow) remain largely independent.
A notable illustration of this is the treatment of power generation cost functions, where nonlinear network optimization has largely used polynomial representations and market optimization has adopted piecewise linear encodings.
This work combines state-of-the-art results from both lines of research to understand the best mathematical formulations of the nonlinear AC optimal power flow problem with piecewise linear generation cost functions.
An extensive numerical analysis of non-convex models, linear approximations, and convex relaxations across fifty-four realistic test cases illustrates that nonlinear optimization methods are surprisingly sensitive to the mathematical formulation of piecewise linear functions.  
The results indicate that a poor formulation choice can slow down algorithm performance by a factor of ten, increasing the runtime from seconds to minutes.
These results provide valuable insights into the best formulations of nonlinear optimal power flow problems with piecewise linear cost functions, an important step towards building a new generation of energy markets that incorporate the nonlinear AC power flow model.
\end{abstract}



\section*{Nomenclature}
\addcontentsline{toc}{section}{Nomenclature}
\begin{IEEEdescription}[\IEEEusemathlabelsep\IEEEsetlabelwidth{$Y^s = g^s-$}]
  \item [{$N$}]  - The set of nodes 
  \item [{$E$, $E^R$}]  - The set of {\em from} and {\em to} branches 
  \item [{$G$}]  - The set of generators
  \item [{$G_i$}] - The subset of generators at bus $i$
  \item [{$p_k$}] - The number of points in piecewise linear cost for generator $k$
  \item [{$C_k$}] - The cost points $[1,2,\dots,p_k]$ for generator $k$
  \item [{$C'_k$}] - The cost points $[2,3,\dots,p_k]$ for generator $k$
  \item [{$C''_k$}] - The cost points $[3,4,\dots,p_k]$ for generator $k$
  \item [{$cg_{kl}$}] - The cost of generator $k$ at point $l \in C_k$ 
  \item [{$pg_{kl}$}] - The power of generator $k$ at point $l \in C_k$
  \item [{${\Delta cg}_{kl}$}] - The incremental cost of generator $k$ between points $l \in C'_k$ and $l-1$
  \item [{$bcg_{kl}$}] - The cost offset of generator $k$ between points $l \in C'_k$ and $l-1$
  \item [{$a_k,b_k,c_k$}] - Polynomial cost coefficients of generator $k$
  \item [{$\bm i$}] - Imaginary number constant
  \item [{$S = p+ \bm iq$}] - AC power
  \item [{$V = v \angle \theta$}]  - AC voltage
  \item [{$Y = g + \bm ib$}]  - Branch admittance
  \item [{$W$}]  - Product of two AC voltages
  \item [{$s^u$}] - Branch apparent power thermal limit
  \item [{$\theta^{\Delta l}, \theta^{\Delta u}$}] - Voltage angle difference limits
  \item [{$S^g$, $S^d$}] - AC power generation and demand
  \item [{$\Re(\cdot), \Im(\cdot)$}] - Real and imaginary parts of a complex num.
  \item [{$(\cdot)^*$, $|\cdot|$}] - Conjugate and magnitude of a complex num.
  \item [{$x^l, x^u$}] - Lower and upper bounds of $x$, respectively
  \item [{$\bm x$}] - A constant value
\end{IEEEdescription}

\section{Introduction}

\IEEEPARstart{O}{ver} the last several decades, competitive energy markets have proven to be an effective mechanism to generate power at minimal cost. 
In such markets potential energy generation units provide offers in the form of piecewise convex functions of generation production cost, and an independent system operator (ISO) or regional transmission organization (RTO) solves a mathematical optimization problem to determine the cheapest dispatch of those generating units, considering a wide variety of network reliability criteria, that is market clearing~\cite{carlson2012miso,wang2013extreme,chen2017mip}.
Given the uncertainty of future energy demands, the market clearing process is repeated at different time scales ranging from day-head to real-time, at approximately 1 hour and 15 minute intervals, respectively.

The significant size of real-world market clearing problems, which feature 100s to 1000s of generating units, and tight run-time requirements, just a few minutes, present a significant computational challenge to optimization algorithms.  
Modern energy markets currently solve this challenging optimization problem by utilizing commercial mixed-integer linear optimization software, such as CPLEX \cite{cplex}, Gurobi \cite{gurobi}, and Xpress \cite{xpress},  
all of which provide high reliability and state-of-the-art computational performance.
However, a key limitation of these commercial tools is a focus on linear equations, limiting market clearing optimization to linear approximations of power flow physics, such as the seminal DC Power Flow \cite{4956966}.
This approximation of the true nonlinear physics of AC power networks results in out-of-market corrections by network operators to adjust for inaccuracies in the market's physics model. 

Recent advances in nonlinear optimization \cite{Ipopt} and convex nonlinear relaxations \cite{8635446} have spurred aspirations for a new generation of market clearing optimization software that considers the full AC power flow physics.
An AC market design has the promise of both reducing out-of-market corrections and incorporating prices for valuable ancillary services such as voltage support capabilities, which are currently priced by ad-hoc methods.
From an algorithmic standpoint, the realization of an AC market in practice requires the fast and reliable solution of challenging mixed-integer non-convex nonlinear optimization problems, which is an active area of research in the optimization community.
Recognizing the near-term potential for AC power flow markets, in 2019 ARPA-e conducted a Grid Optimization Competition \cite{goc}, to identify the most promising algorithmic approaches for building the next generation of power network optimization software.

Pursuing a future AC market design, this work explores how to best model piecewise convex functions of generation production costs in nonlinear optimization algorithms.
Specifically, it considers different mathematically-equivalent formulations of piecewise linear cost functions and evaluates their computational implications on nonlinear power network optimization.
The core observations of this work are threefold: (1) the mathematical modeling lessons learned from linear active-power-only markets do not necessarily carry over to forthcoming AC power markets; (2) a poor choice of the piecewise linear cost function representation can result in a solution time slowdown of as much as 10 times on realistic test systems; (3) ultimately, the ``$\lambda$'' and ``$\Delta$'' formulations of the piecewise cost functions prove to be the most suitable for current nonlinear optimization software. 
To the best of our knowledge, this paper provides the first synthesis of state-of-the-art models for piecewise linear cost functions coming from the market operations literature \cite{knueven2020mixed} with state-of-the-art models for nonlinear AC optimal power flow literature \cite{7271127,pm_pscc}.

The rest of the paper is organized as follows: 
Section \ref{sec:opf} provides a brief introduction to optimal power flow and Section \ref{sec:pwl} reviews mathematical models for piecewise linear cost functions.  Section \ref{sec:computation} conducts computational experiments of 12 optimization models across 54 test cases providing the core contributions of the paper, and Section \ref{sec:conclusion} finishes with closing remarks.

\section{Optimal Power Flow}
\label{sec:opf}

The most prevalent power network optimization task is arguably the optimal power flow problem (OPF).  At a high level, the OPF problem is a single-time period optimization task that consists of finding the cheapest way to generate sufficient power to meet a specified demand.  The challenge of OPF is that an AC power network with a variety of operating constraints is used to transmit the power from the generators to the demands.  Capturing both these operational constraints and the AC physics of a power network gives rise to a challenging non-convex nonlinear optimization task and a thriving body of algorithmic research \cite{8635446}.  It is important to highlight, although OPF forms a foundational sub-problem of an energy market, that a real-world market requires many extensions including co-optimization of multiple-time periods, commitment of generation units, and contingency constraints, to name a few.  This work focuses on OPF as a first necessary step to building a comprehensive and reliable AC energy market optimization.

This work begins with the canonical academic AC-OPF formulation \cite{1908.02788} presented in Model \ref{model:ac_opf}.  At a high level, the network is defined by a set of buses $N$, lines $E$, and generators $G$.  A notable feature of this formulation is that the active power generator costs are provided as convex quadratic functions, that is,
\begin{align}
\mbox{Cost}_k(x) = \bm {c_k} x^2 + \bm {b_k} x + \bm {a_k}  \;\; \forall k \in G \label{eq:cost_poly}
\end{align}
The objective function \eqref{ac_opf_poly_obj} seeks to minimize total of the generation costs.  
The constraints \eqref{ac_opf_1},\eqref{ac_opf_2} capture the bus voltage requirements and generator output limits, respectively.  The constraints \eqref{ac_opf_3},\eqref{ac_opf_4} model the AC power flow physics via power balance and Ohm's law, respectively.  Finally, the constraints \eqref{ac_opf_5},\eqref{ac_opf_6} enforce the thermal and angle stability limits across the power lines in the network.  For additional details, a first-principles derivation of this model is presented in \cite{7271127}.

\begin{model}[t]
\caption{AC Optimal Power Flow (AC-OPF)}
\label{model:ac_opf}
\begin{subequations}
\vspace{-0.2cm}
\begin{align}
& \mbox{\bf variables: } S^g_k (\forall k \in G), V_i ( \forall i \in N), S_{ij} (\forall (i,j) \in E \cup E^R) \nonumber \\
%
& \mbox{\bf minimize:} \sum_{k \in G} \mbox{Cost}_k(\Re(S^g_k)) \label{ac_opf_poly_obj} \\
&\mbox{\bf subject to:} \nonumber \\
& \bm {v^l}_i \leq |V_i| \leq \bm {v^u}_i \;\; \forall i \in N \label{ac_opf_1} \\
& \bm {S^{gl}}_k \leq S^g_k \leq \bm {S^{gu}}_k \;\; \forall k \in G \label{ac_opf_2}  \\
& \sum_{k \in G_i} S^g_k - \bm S^d_i = \!\!\!\!\!\!\!\!\!\! \sum_{\substack{(i,j)\in E_{i} \cup E^R_{i}}} \!\!\!\!\!\!\!\!\! S_{ij} \; \forall i\in N \label{ac_opf_3} \\ 
& S_{ij} = \bm Y^*_{ij} |V_i|^2 - \bm Y^*_{ij} V_i V_j^* \;\; (i,j)\in E \cup E^R \label{ac_opf_4}\\
& |S_{ij}| \leq \bm {s^u}_{ij} \;\; \forall (i,j) \in E \cup E^R \label{ac_opf_5}  \\
& \bm {\theta^{\Delta l}}_{ij} \leq \angle (V_i V^*_j) \leq \bm {\theta^{\Delta u}}_{ij} \;\; \forall (i,j) \in E \label{ac_opf_6}
\end{align}
\end{subequations}
\end{model}

\subsection{OPF Solution Methods}

The mathematical optimization problem presented by Model \ref{model:ac_opf} is a non-convex NonLinear Program (NLP).
This problem is known to feature local minima \cite{6581918} and to be NP-Hard in the general case \cite{verma2009power,7063278}.
The standard solution for challenging non-convex NLP models is use of global optimization solvers, such as {\sc Baron} \cite{ts:05}, Couenne \cite{Belotti:2009}, and Alpine \cite{NagarajanLuWangBentSundar2019}, which provide solution quality bounds and optimality proofs.
However, these approaches can only solve AC-OPF problems with a few hundred buses, which is an order of magnitude less than real-world applications that require thousands of buses.
It has been observed that interior point methods, such as Ipopt \cite{Ipopt} and {\sc Knitro} \cite{Byrd2006}, can quickly find high-quality solutions to real-world AC-OPF problems \cite{7271127}.
However, these methods do not provide global guarantees of convergence and optimality of the solutions that they find.
This performance gap between off-the-self interior point and global optimization solvers has yielded  a wide variety of bespoke solutions for solving AC-OPF with quality guarantees, including problem-specific convex relaxations \cite{1664986,7271127}, polynomial optimization \cite{7038397}, and bound tightening \cite{cp_qc_fp}.
See \cite{8635446} for a comprehensive review of different approaches.

In addition to the AC-OPF formulation presented in Model \ref{model:ac_opf}, this work also considers two canonical alternatives:
the seminal DC Power Flow approximation and a simple convex relaxation of the OPF problem.
These alternatives serve to position the results of this work in the ongoing transition of power network optimization from linear active-power-only approximations to nonlinear active-and-reactive optimization.

\subsection{Convex Relaxation}
Convex relaxations of the non-convex AC model have drawn significant interest in recent years \cite{8635446}, in large part due to their ability to provide tight bounds on the AC-OPF solution quality.  Following those lines, this work considers the Second-Order Cone (SOC) relaxation of the AC power flow equations \cite{1664986}.  Although some relaxations are stronger than SOC \cite{7271127,6345272,moment_hierarchy} and others are faster than SOC \cite{7540869}, the SOC relaxation is selected because it provides an appealing tradeoff between bounding strength and runtime performance.

The first insight of the SOC relaxation is that the voltage product terms $V_i^* V_j$ can be lifted into a higher dimensional $W$-space as follows:
\begin{subequations}
\begin{align}
|V_i|^2 &\Rightarrow W_{ii} \;\; \forall i \in N \label{w_lift_1} \\
V_i V^*_j &\Rightarrow W_{ij} \;\; \forall (i,j) \in E \label{w_lift_2}
\end{align}
\end{subequations}
Note that lifting Model \ref{model:ac_opf} into the $W$-space makes all of the non-convex constraints linear.  
The second insight of the SOC relaxation is that the $W$-space relaxation can be strengthened by adding the valid inequality,
%
\begin{align}
& |W_{ij}|^2 \leq W_{ii}W_{jj} \;\; \forall (i,j) \in E \label{w_soc}
\end{align}
%
which is a convex rotated SOC constraint that is supported by a wide variety of commercial optimization tools.

Utilizing these two insights, the SOC-OPF relation of the AC-OPF problem is presented in Model \ref{model:soc_opf}.  Many of the constraints remain the same; the core differences are as follows.  Constraints \eqref{soc_opf_1},\eqref{soc_opf_2},\eqref{soc_opf_3},\eqref{soc_opf_4} capture the bus voltage requirements, Ohm's law, and voltage angle stability limits in the lifted $W$-space.  Constraint \eqref{soc_opf_5} is a new constraint that strengthens the relaxation.
The virtues of this model are that it is convex (i.e., global optimality is achieved by local solvers like Ipopt and {\sc Knitro}) and it provides a lower bound to the objective function value of the non-convex AC-OPF model.  The principal weakness of this model is that its voltage solution is non-physical and usually does not provide useful insights into the solution of the non-convex AC-OPF model.

\begin{model}[t]
\caption{SOC Optimal Power Flow (SOC-OPF)}
\label{model:soc_opf}
\begin{subequations}
\vspace{-0.2cm}
\begin{align}
& \mbox{\bf variables: } S^g_k (\forall k \in G), S_{ij} (\forall (i,j) \in E \cup E^R), \phantom{123456789} \nonumber \\
& W_{ii} (\forall i \in N), W_{ij} (\forall i,j \in E) \nonumber \\
& \mbox{\bf minimize: } \eqref{ac_opf_poly_obj} \nonumber \\
&\mbox{\bf subject to: } \eqref{ac_opf_2}, \eqref{ac_opf_3}, \eqref{ac_opf_5} \nonumber \\
& (\bm {v^l}_i)^2 \leq W_{ii} \leq (\bm {v^u}_i)^2 \;\; \forall i \in N \label{soc_opf_1} \\
& S_{ij} = \bm Y^*_{ij} W_{ii} - \bm Y^*_{ij} W_{ij} \;\; (i,j)\in E \label{soc_opf_2}\\
& S_{ji} = \bm Y^*_{ij} W_{jj} - \bm Y^*_{ij} W_{ij}^* \;\; (i,j)\in E \label{soc_opf_3}\\
& \tan(\bm {\theta^{\Delta l}}_{ij})\Re(W_{ij}) \leq \Im(W_{ij}) \nonumber \\ & \leq \tan(\bm {\theta^{\Delta u}}_{ij})\Re(W_{ij}) \;\; \forall (i,j) \in E \label{soc_opf_4} \\
& |W_{ij}|^2 \leq W_{ii}W_{jj} \;\; (i,j)\in E \label{soc_opf_5}
\end{align}
\end{subequations}
\end{model}

\begin{model}[t]
\caption{DC Optimal Power Flow (DC-OPF)}
\label{model:dc_opf}
\begin{subequations}
\vspace{-0.2cm}
\begin{align}
& \mbox{\bf variables: } p^g_k (\forall k \in G), p_{ij} (\forall (i,j) \in E \cup E^R), \theta_i (\forall i \in N) \nonumber \\
%
& \mbox{\bf minimize:} \sum_{k \in G} \mbox{Cost}_k(p^g_k) \label{dc_opf_poly_obj} \\
& \mbox{\bf subject to:} \nonumber \\
& \Re(\bm {S^{gl}}_k) \leq p^g_k \leq \Re(\bm {S^{gu}}_k) \;\; \forall k \in G \label{dc_opf_1}  \\
& \sum_{k \in G_i} p^g_k - \Re(\bm S^d_i) = \!\!\!\!\!\!\!\!\!\! \sum_{\substack{(i,j)\in E_{i} \cup E^R_{i}}} \!\!\!\!\!\!\!\!\! p_{ij} \; \forall i\in N \label{dc_opf_2} \\ 
& p_{ij} = - \Im(\bm Y^*_{ij})(\theta_i - \theta_j) \;\; (i,j)\in E \cup E^R \label{dc_opf_3}\\
& -\bm {s^u}_{ij} \leq p_{ij} \leq \bm {s^u}_{ij} \;\; \forall (i,j) \in E \cup E^R \label{dc_opf_4}  \\
& \bm {\theta^{\Delta l}}_{ij} \leq \theta_i - \theta_j \leq \bm {\theta^{\Delta u}}_{ij} \;\; \forall (i,j) \in E \label{dc_opf_5}
\end{align}
\end{subequations}
\end{model}

\subsection{Linear Approximation}

The most widely used approach to solving OPF problems is to approximate the power flow physics with the DC Power Flow approximation \cite{4956966}.  This model is achieved by taking a first-order Taylor expansion around the nominal voltage operating point of $V_i \approx 1 \; \angle \; 0$.  This expansion yields the following voltage product approximation,
\begin{align}
V_i V^*_j &\Rightarrow 1 \; \angle \; (\theta_i - \theta_j)
\end{align}
and results in omitting the reactive power variables and constraints from the model as they are constant values in the first-order Taylor expansion.
A more detailed derivation of this model from first principles is available in \cite{4956966,LPAC_ijoc}.

Applying this transformation to the AC-OPF problem yields the DC-OPF approximation in Model \ref{model:dc_opf}.  
The constraints are similar to the AC-OPF model but with the reactive power components omitted.
The objective function \eqref{dc_opf_poly_obj} is the same but applied to the active-power variables directly.
The constraint \eqref{dc_opf_1} captures the generator output limits.
The constraints \eqref{dc_opf_2},\eqref{dc_opf_3} approximate the power flow physics of power balance and Ohm's law, respectively.
Finally, the constraints \eqref{dc_opf_4},\eqref{dc_opf_5} enforce the thermal and angle stability limits across the power lines in the network.
The virtue of this model is that it is a linear optimization problem, which benefits from decades of research and mature commercial optimization software.  The principal weakness of this model is that it cannot consider bus voltage and reactive power requirements.

\section{Piecewise Linear Costs}
\label{sec:pwl}

Mathematical formulations for piecewise linear cost functions are a cornerstone of linear and mixed-integer linear optimization tools \cite{bertsimas1997introduction,wolsey1999integer}.
The first linear programming formulation for convex piecewise linear costs appeared shortly after the simplex method was popularized for a practical reason -- to approximate separable convex functions \cite{charnes1954minimization}. 
This was quickly followed by several alternative formulations for separable convex piecewise linear costs \cite{dantzig1956recent,dantzig1958linear,ho1985relationships}.
Formulations for non-convex piecewise linear cost functions arrived shortly after the development of Mixed-integer linear programming (MILP), a much richer modeling framework, that can capture many discrete mathematical structures \cite{dantzig1960significance}.
Both the convex and non-convex formulations of piecewise linear cost have been a fruitful research topic over the years.  See \cite{vielma2010mixed} for a recent comprehensive review of different approaches.

In the context of power systems, the majority of the piecewise linear cost literature has focused on the formulation of convex piecewise linear cost functions in the context of the unit commitment (UC) problem, which focuses on the temporal constraints of assigning generators to deliver power for several hours or days. 
The use of convex piecewise linear cost in the UC context dates back to Garver \cite{garver1962power}, who proposed them in a MILP formulation for UC. 
Over the years research demonstrated that the formulation of the piecewise linear production cost can have a profound impact on modern MILP solver performance \cite{sridhar2013locally,chen2017mip,knueven2020mixed}.
While several papers have studied quadratic and polynomial convex production costs for UC \cite{frangioni2006solving,frangioni2008tighter,bacci2019new}, these formulations have not been adopted by industry in the United States, mainly due to the increased complexity of the resulting nonlinear mixed-integer optimization problem.
At this time, it is still the standard practice in the United States for ISOs and RTOs to require generators to submit convex piecewise linear offers for generation production costs.

As previously discussed, there is an increasing interest for market operators to consider AC physics directly in market clearing problems, like OPF and UC \cite{goc}.
However, there is currently an inconsistency between the existing market structures, which model generator costs as convex piecewise linear functions, and the AC-OPF literature that has standardized around convex quadratic cost functions \cite{1908.02788}.
The principal objective of this work is to consider the breadth of convex piecewise linear cost formulations developed in the UC literature in the context of AC-OPF, to understand their performance implications on forthcoming nonlinear optimization problems.

\begin{figure}[t]
    \begin{center}
    \includegraphics[width=0.48\textwidth]{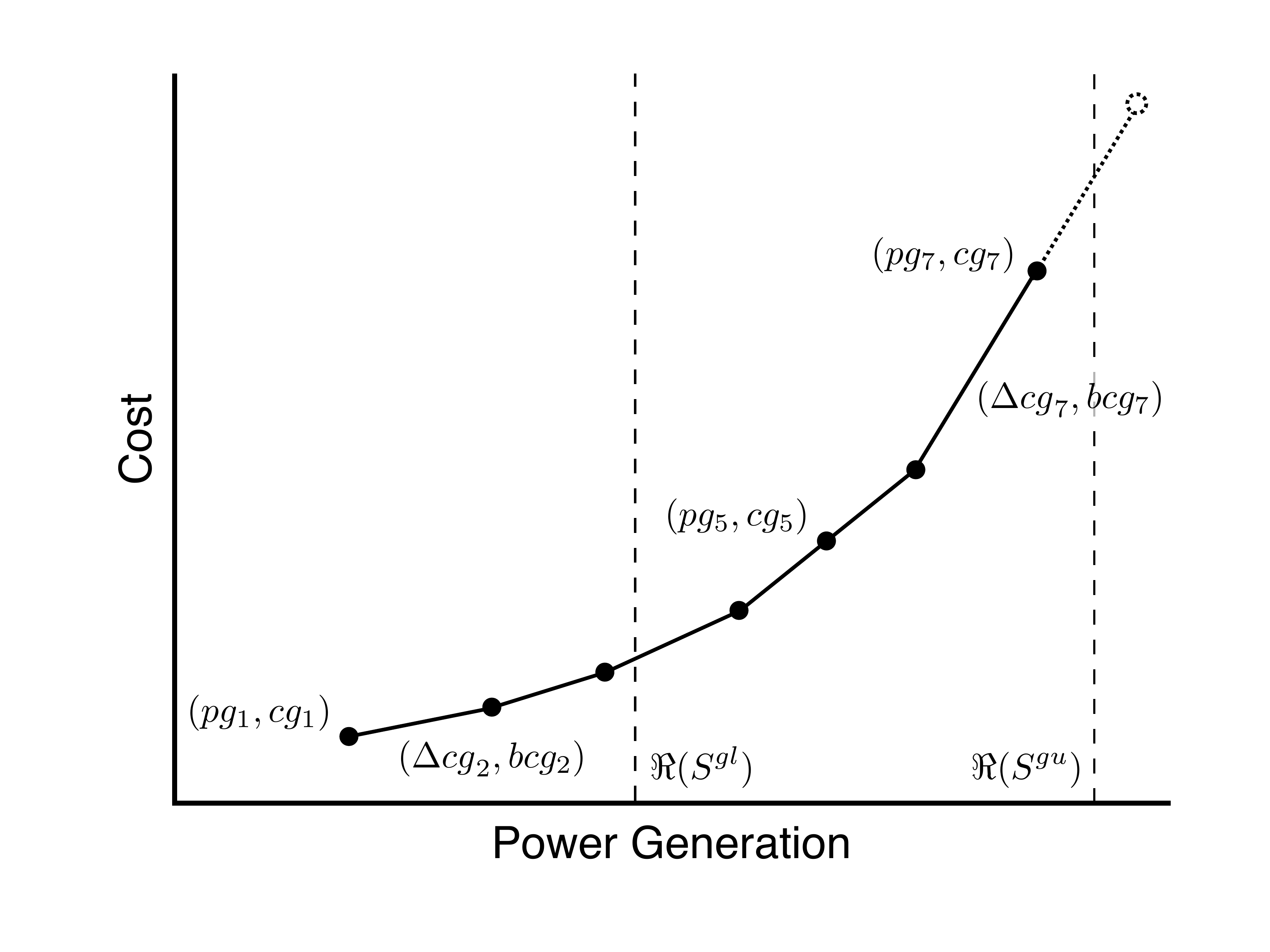}
    \end{center}
    \vspace{-0.3cm}
    \caption{An example of typical piecewise linear data processing.  In this case, the points 1 and 2 can be removed due to the generation lower bound, $\Re(S^{gl})$.  Point 5 can be removed due to no change in the slope of the adjacent segments, and point 7 needs to be extended to include the generation upper bound, $\Re(S^{gu})$.}
    \label{fig:pwl_points}
\end{figure}

\subsection{Piecewise Linear Data Model and Assumptions}

In this work piecewise linear functions are defined by a sequence of points representing line segments for each generator $k \in G$,
%
\begin{align}
(pg_{kl}, cg_{kl}) \;\; \forall l \in C_k \label{eq:pwl_pt}
\end{align}
%
where $cg_{kl}$ is the cost of generating $pg_{kl}$ megawatts of power and the set $C_k$ determines how many points each generator cost function has.
In some cases it is convenient to consider the piecewise linear function as a collection of lines defined by slope-intercept pairs as follows,
\begin{subequations}
\begin{align}
{\Delta cg}_{kl} = \frac{cg_{k,l} - cg_{k,l-1}}{pg_{k,l} - pg_{k,l-1}} \;\; \forall l \in C'_k \label{eq:pwl_si_1} \\
bcg_{kl} = cg_{kl} - {\Delta cg}_{kl}pg_{kl}   \;\; \forall l \in C'_k \label{eq:pwl_si_2}
\end{align}
\end{subequations}
where ${\Delta cg}_{kl}$ and $bcg_{kl}$ are the slope and intercept respectively.  Note that the set $C'_k$ omits the first index of the standard point set $C_k$ so that $l-1$ is well defined.
Figure \ref{fig:pwl_points} provides an illustration of a prototypical generation cost function encoded as a sequence of points.

{\em Assumptions:}
A core challenge of working with piecewise linear functions is their generality.
Hence, it is important to provide a detailed specification of the inputs that are permitted.
Throughout this work it is assumed that the optimization problem of interest is a minimization problem and the points of each generator encode a convex piecewise linear cost function.  Specifically, this work requires three key properties of these functions: (1) the generator's lower bound, $S^{gl}_k$, occurs in the first segment of function; (2) the generator's upper bound, $S^{gu}_k$, occurs in the last segment of function; (3) the slope of each linear segment is strictly increasing.  These properties are summarized as follows,
\begin{subequations}
\begin{align}
pg_{k,1} & \leq \Re(S^{gl}_k) < pg_{k,2} & \;\; \forall k \in G \label{eq:pwl_as_1} \\
pg_{k,{p_k}-1} & < \Re(S^{gu}_k) \leq pg_{k,p_k} & \;\; \forall k \in G \label{eq:pwl_as_2} \\
{\Delta cg}_{k,l-1} & < {\Delta cg}_{k,l} &  \forall l \in C''_k \;\; \forall k \in G \label{eq:pwl_as_3}
\end{align}
\end{subequations}
If the given piecewise linear function is convex, these properties can be ensured by the following data processing procedure.

{\em Data Processing:}
For a variety of reasons, real-world piecewise linear generation cost functions often benefit from data cleaning before encoding them in an optimization model.  In this work we conduct the following  data processing procedures: (1) if the generator bounds are outside of the first and last  segments, that is $\Re(S^{gl}_k) < pg_{k,1}$ or $pg_{k,p_k} < \Re(S^{gu}_k)$, then the segments are extended to include the generator bounds; (2) if the piecewise linear function includes segments beyond the bounds of the generator, that is $pg_{k,2} \leq \Re(S^{gl}_k)$ or $\Re(S^{gu}_k) \leq pg_{k,{p_k}-1}$, then the extra out-of-bounds segments are removed; (3) if there is little change in the slope of two adjacent segments, that is ${\Delta cg}_{k,l-1} \approx {\Delta cg}_{k,l}$, then they are combined into one segment.  These simple data processing steps increase the performance of optimization algorithms by removing redundant constraints and serve to enforce the mathematical requirements of the formulations considered by this work.
Figure \ref{fig:pwl_points} provides an illustration of generation cost function that requires these data processing procedures.

\subsection{Formulations of OPF with Piecewise Linear Costs}

The interest of this work is variants of the OPF problems from Section \ref{sec:opf} where the polynomial active power generation cost functions are replaced with piecewise linear functions, that is,
\begin{align}
\mbox{Cost}_k(x) = \max_{l \in C'_k} \{\bm {\Delta cg}_{kl} x + \bm {bcg}_{kl}\}  \;\; \forall k \in G \label{eq:cost_pwl}
\end{align}
Unlike the previously considered polynomial functions, there are a wide variety of equivalent mathematical encodings of piecewise linear functions.
Following the seminal works on piecewise linear formulations \cite{ho1985relationships,fourer1992simplex}, this work considers the four ``standard'' formulations for convex piecewise linear functions, referred to as the $\Psi$, $\lambda$, $\Delta$ and $\Phi$ models by \cite{fourer1992simplex}.  All four formulations are mathematically equivalent but can have significant performance implications for the numerical methods used in optimization algorithms.
While other formulations for piecewise linear costs exist, three of these standard formulations are commonly used in power systems problems (i.e., $\Psi$, $\lambda$, $\Delta$) \cite{knueven2020mixed}. This work includes a fourth formulation, $\Phi$, for completeness and because it shares an interesting mathematical connection as the dual of the $\Delta$ formulation.  
At a high level, these formulations represent two distinct perspectives on modeling piecewise linear functions, focusing on either the function evaluation points (i.e., $\Psi$,$\lambda$) or integration of the cost function's derivative (i.e., $\Delta$,$\Phi$).

{\em The $\Psi$ Formulation:}
This formulation is arguably the most popular and intuitive. The $\Psi$ formulation explicitly models the epigraph of \eqref{eq:cost_pwl}, that is the region on and above the objective function on the graph ($x$, $f(x)$).
However, because mathematical programming solvers do not usually have explicit support for the $\max$ function, an auxiliary cost variable $c^g_k \in [\bm {cg}_{k,1}, \bm {cg}_{k,p_k}]$ is introduced for each generator, which combined with inequality constraints, captures the semantics of the $\max$ function.
The complete AC-OPF problem using the $\Psi$ cost formulation is presented in Model \ref{model:ac_opf_psi}.

\begin{model}[t]
\caption{AC-OPF with the $\Psi$ Cost Model (AC-OPF-$\Psi$)}
\label{model:ac_opf_psi}
\begin{subequations}
\vspace{-0.2cm}
\begin{align}
& \mbox{\bf variables: } S^g_k (\forall k \in G), S_{ij} (\forall (i,j) \in E \cup E^R), V_i (\forall i \in N), \nonumber \\
& c^g_k \in [\bm {cg}_{k,1}, \bm {cg}_{k,p_k}] \; (\forall k \in G) \nonumber \\
& \mbox{\bf minimize: } \sum_{k \in G} c^g_k \label{ac_opf_psi_obj} \\
& \mbox{\bf subject to: } \eqref{ac_opf_1}-\eqref{ac_opf_6} \nonumber \\
& c^g_k \geq \bm {\Delta cg}_{kl} \Re(S^g_k) + \bm {bcg}_{kl} \;\; \forall l \in C'_k, \forall k \in G \label{ac_opf_psi_1} 
\end{align}
\end{subequations}
\end{model}

{\em The $\lambda$ Formulation:}
This formulation reflects the most natural encoding of a convex hull from the collection of points \eqref{eq:pwl_pt}, which is a popular modeling approach in the linear programming literature \cite{dantzig1998linear}.
The core idea is to introduce an interpolation variable $\lambda^{cg}_{kl} \in [0, 1]$ for each point in the piecewise linear function and link all of these interpolation variables together with the constraint $\sum_{l \in C_k} \lambda^{cg}_{kl} = 1$.
The power and cost of the interpolated point can be recovered with the expressions $\sum_{l \in C_k} \bm {pg}_{kl} \lambda^{cg}_{kl}$ and $\sum_{l \in C_k} \bm {cg}_{kl} \lambda^{cg}_{kl}$, respectively.
The complete AC-OPF problem using the $\lambda$ cost formulation is presented in Model \ref{model:ac_opf_lam}.  
The $\lambda$ formulation is also interesting because it represents the mathematical dual of the $\Psi$ piecewise linear formulation \cite{fourer1992simplex}.

\begin{model}[t]
\caption{AC-OPF with the $\lambda$ Cost Model (AC-OPF-$\lambda$)}
\label{model:ac_opf_lam}
\begin{subequations}
\vspace{-0.2cm}
\begin{align}
& \mbox{\bf variables: } S^g_k (\forall k \in G), S_{ij} (\forall (i,j) \in E \cup E^R), V_i (\forall i \in N), \nonumber \\
& \lambda^{cg}_{kl} \in [0, 1] \; (\forall l \in C_k, \forall k \in G) \nonumber \\
&\mbox{\bf minimize: } \sum_{k \in G} \sum_{l \in C_k} \bm {cg}_{kl} \lambda^{cg}_{kl} \label{ac_opf_lam_obj} \\
& \mbox{\bf subject to: }  \eqref{ac_opf_1}-\eqref{ac_opf_6} \nonumber \\
& \sum_{l \in C_k} \bm {pg}_{kl} \lambda^{cg}_{kl} = \Re(S^g_k) \;\; \forall k \in G \label{ac_opf_lam_1} \\
& \sum_{l \in C_k} \lambda^{cg}_{kl} = 1 \;\; \forall k \in G \label{ac_opf_lam_2}
\end{align}
\end{subequations}
\end{model}

{\em The $\Delta$ Formulation:}
This formulation breaks the cost function into a collection of generation bins, ${\Delta pg}_{kl}$, based on consecutive points in the piecewise linear function, that is, ${\Delta pg}_{kl} \in [0, \bm {pg}_{kl} - \bm {pg}_{k,l-1}]$.
The key observation is that each bin can be associated with a linear cost based on the slope of that line segment, i.e., $\bm {\Delta cg}_{kl} {\Delta pg}_{kl}$.
An interpretation of this formulation is that it computes the integral of the cost's derivative along the power-axis of the piecewise linear function.
The complete AC-OPF problem using the $\Delta$ cost formulation is presented in Model \ref{model:ac_opf_delta}.
Both the $\Delta$ formulation and the $\lambda$ formulation are particularly interesting as they have proven the most effective formulations in practical UC problems \cite{knueven2020mixed}. 

\begin{model}[t]
\caption{AC-OPF with the $\Delta$ Cost Model (AC-OPF-$\Delta$)}
\label{model:ac_opf_delta}
\begin{subequations}
\vspace{-0.2cm}
\begin{align}
& \mbox{\bf variables: } S^g_k (\forall k \in G), S_{ij} (\forall (i,j) \in E \cup E^R), V_i (\forall i\in N), \nonumber \\
& {\Delta pg}_{kl} \in [0, \bm {pg}_{k,l} - \bm {pg}_{k,l-1}] \; (\forall l \in C'_k, \forall k \in G) \nonumber \\
& \mbox{\bf minimize: } \sum_{k \in G} \left( \bm {cg}_{k,1} + \sum_{l \in C'_k} \bm {\Delta cg}_{kl} {\Delta pg}_{kl} \right) \label{ac_opf_delta_obj} \\
& \mbox{\bf subject to: } \eqref{ac_opf_1}-\eqref{ac_opf_6} \nonumber \\
& \Re(S^g_k) = \sum_{l \in C'_k} {\Delta pg}_{kl} + \bm {pg}_{k,1} \;\; \forall k \in G \label{ac_opf_delta_1} 
\end{align}
\end{subequations}
\end{model}

{\em The $\Phi$ Formulation:}
This formulation is the most challenging to interpret and combines elements from both the $\Psi$ and $\Delta$ models.
It begins by extracting a linear cost function from the first segment of the piecewise linear function, i.e., $\bm {\Delta cg}_{k,2}\Re(S^g_k) + \bm {bcg}_{k,2}$.  It then defines bins for how much power is supplied by each segment after the first, i.e., $\Phi_{kl} \in [0, \Re (\bm {S^{gu}}_k) - \bm {pg}_{k,l-1}]$.  The extra incremental cost of each segment over the previous segments is then captured by $(\bm {\Delta cg}_{k,l} - \bm {\Delta cg}_{k,l-1}) \Phi_{kl}$.
An interpretation of this formulation is that it computes the integral of the cost's derivative along the cost-axis of the piecewise linear function.
The complete AC-OPF problem using the $\Phi$ cost formulation is presented in Model \ref{model:ac_opf_phi}.
Although this $\Phi$ formulation is uncommon in the literature, it is interesting in this context as the mathematical dual of the $\Delta$ piecewise linear formulation \cite{fourer1992simplex}.

\begin{model}[t]
\caption{AC-OPF with the $\Phi$ Cost Model (AC-OPF-$\Phi$)}
\label{model:ac_opf_phi}
\begin{subequations}
\vspace{-0.2cm}
\begin{align}
& \mbox{\bf variables: } S^g_k (\forall k \in G), S_{ij} (\forall (i,j) \in E \cup E^R), V_i (\forall i\in N \nonumber \\
& \Phi_{kl} \in [0, \Re (\bm {S^{gu}}_k) - \bm {pg}_{k,l-1}]  \;\; \forall l \in C''_k, \forall k \in G \nonumber \\
& \mbox{\bf minimize: }\sum_{k \in G} \left( \bm {\Delta cg}_{k,2}\Re(S^g_k) + \bm {bcg}_{k,2} +  \phantom{\sum_{l \in C''_k} } \right. \nonumber \\ 
& \left. \phantom{123456789---} \sum_{l \in C''_k}  (\bm {\Delta cg}_{k,l} - \bm {\Delta cg}_{k,l-1}) \Phi_{kl} \right) \label{ac_opf_phi_obj} \\
%
& \mbox{\bf subject to: } \eqref{ac_opf_1}-\eqref{ac_opf_6} \nonumber \\
& \Phi_{kl} \geq \Re(S^g_k) - \bm {pg}_{k,l-1} \;\; \forall l \in C''_k, \forall k \in G
\label{ac_opf_phi_1} 
\end{align}
\end{subequations}
\end{model}

\subsection{Relaxation and Approximation Variants}

Models \ref{model:ac_opf_psi} through \ref{model:ac_opf_phi} present variants of the AC-OPF model with different formulations of piecewise linear cost functions.
As these formulations only changed the objective function in that model, it is clear how similar modifications could also be applied to the DC approximation and SOC relaxation models that were presented in Section \ref{sec:opf}.

\section{Computational Evaluation}
\label{sec:computation}

Combining both the power flow formulations from Section \ref{sec:opf} with the piecewise linear formulations from Section \ref{sec:pwl}, this section conducts a detailed computational evaluation of 12 different OPF formulations, i.e., $\{\mbox{AC}, \mbox{SOC}, \mbox{DC}\} \times \{ \Psi, \lambda, \Delta, \Phi \}$.  The overarching observation is that the choice of  piecewise linear formulation has a much more dramatic impact on state-of-the-art nonlinear optimization algorithms than it does for linear optimization algorithms.  This result is demonstrated by three computational studies, the first focusing on a solution quality comparison, the second highlighting the runtime trends of different piecewise linear formulations, and the third comparing different algorithms for solving convex nonlinear optimization problems, which is particularly relevant when considering convex relaxations of OPF problems.


\subsection{Test Cases and Computational Setting}

The traditional AC-OPF benchmark problems, such as {\sc PGLib}-OPF \cite{1908.02788}, have standardized around convex quadratic cost functions, which precludes their use in this work.
To the best of our knowledge, ARPA-e's Grid Optimization Competition Challenge 1 datasets \cite{goc} represent the first comprehensive source of AC-OPF test cases that feature piecewise linear cost functions and hence these were leveraged for building a test suite in this work.
The publicly available {\em Challenge 1 Final Event} data archive consists of 340 AC-OPF cases spanning 18 distinct networks, which range from 500 to 30,000 buses.
This work down selected that collection to a set of 54 representative AC-OPF cases, three representatives from each distinct network.
Finally, the cases were converted from the Grid Optimization Competition data format into the {\sc Matpower} data format \cite{matpower} to be compatible with established state-of-the-art OPF evaluation tools \cite{pm_pscc}.
The specific details of which scenarios were selected are available in Table \ref{tbl:deltas_time_ac}.

The proposed OPF formulations were implemented in Julia v1.5 as an extension to the PowerModels v0.17 \cite{power_models,pm_pscc} framework, which utilizes JuMP v0.21 \cite{DunningHuchetteLubin2017} as a general-purpose mathematical optimization modeling layer.
The NLP formulations were primarily solved with Ipopt \cite{Ipopt} using the HSL MA27 linear algebra solver \cite{hsl_lib} and a cross validation is conducted with {\sc Knitro} v12.2. 
The LP and QCQP formulations were solved with Gurobi v9.0 \cite{gurobi}.
All of the solvers were configured to terminate once the optimality gap was less than $10^{-6}$ without an explicit time limit.  
The evaluation was conducted on HPE ProLiant XL170r servers with two Intel 2.10 GHz CPUs and 128 GB of memory; however, for consistency the algorithms were configured to only utilize one thread.

\subsection{Solution Quality Validation}

The first observation of this experiment is that all of the formulations considered found solutions of identical quality, up to the numerical tolerances of the optimization algorithms.
For the convex DC-OPF and SOC-OPF formulations this serves as an important validation of the implementation's correctness.
Since these problems are convex and all of the piecewise linear formulation are proven to be mathematically equivalent, all models should converge to a consistent globally optimal solution, up to the accuracy of floating point arithmetic.  Detailed results for these models are omitted in the interest of brevity.

A more surprising result is that all four variants of the non-convex AC-OPF formulation also converge to solutions of identical quality, up to the numerical tolerance of the optimization algorithm.  This result is demonstrated in the detailed results presented in Table \ref{tbl:deltas_time_ac}.  The first AC-OPF-$\lambda$ column shows the locally optimal objective value and the following three columns indicate the absolute difference in the objective value from the other three formulations, all of which are well below the optimality tolerance of $10^{-6}$.
These are encouraging results as they suggest the choice of piecewise linear formulation in nonlinear optimization algorithms can be taken solely on the criteria of runtime performance without concern for solution quality degradation.

\begin{figure}[!ht]
    \begin{center}
    \includegraphics[width=0.50\textwidth]{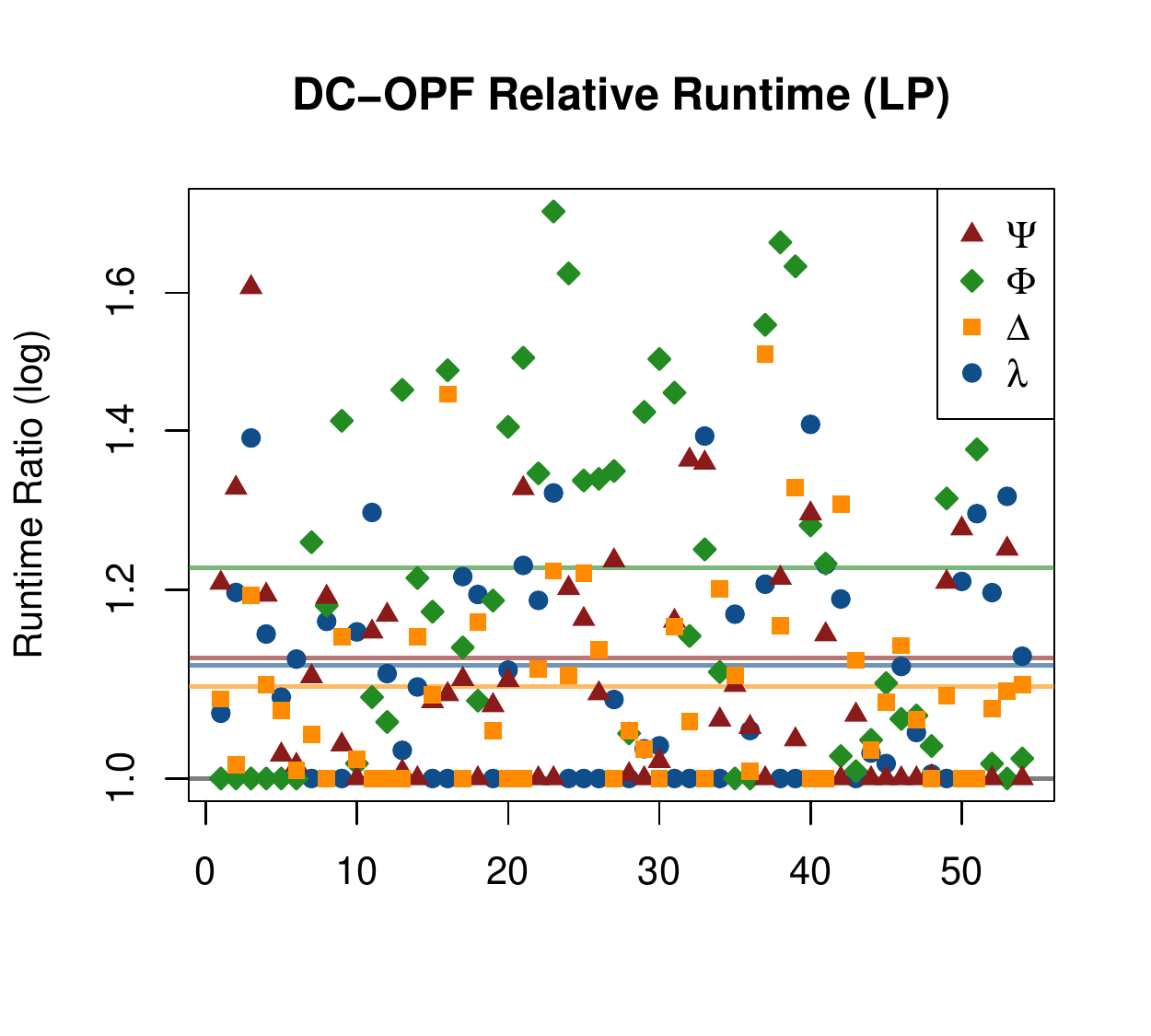}\\
    \vspace{-1.0cm}
    \includegraphics[width=0.50\textwidth]{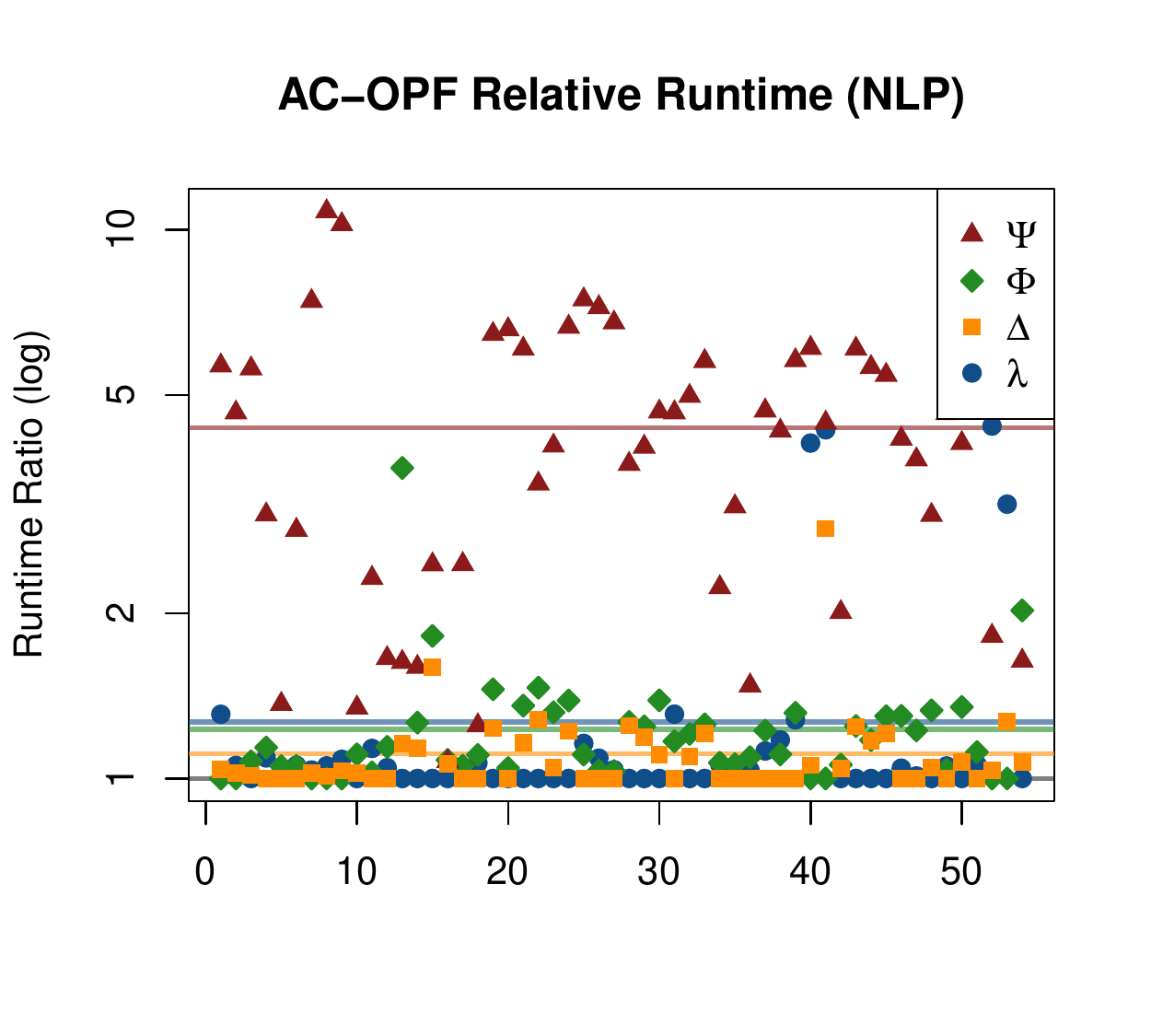}\\
    \vspace{-1.0cm}
    \includegraphics[width=0.50\textwidth]{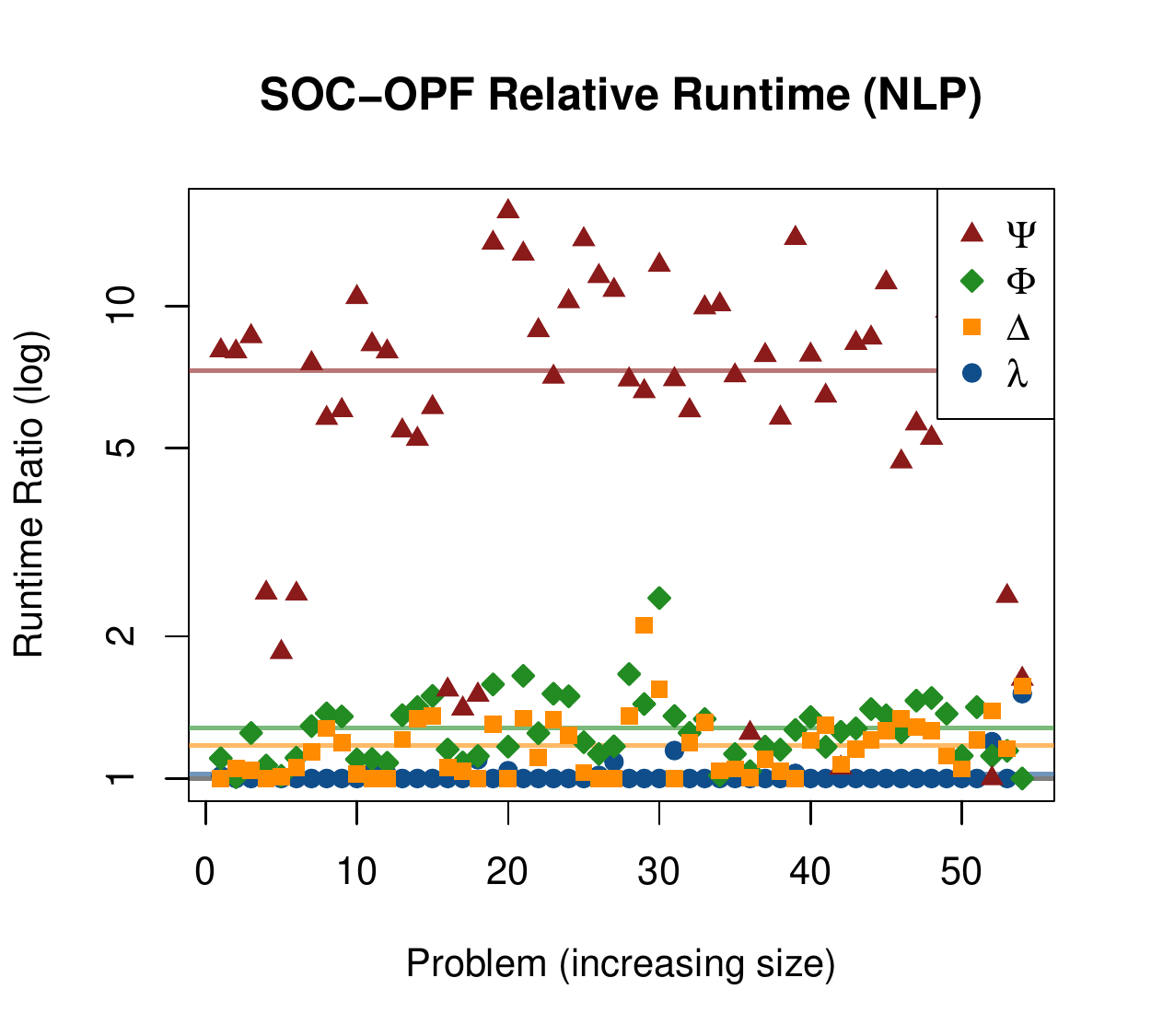}
    \end{center}
    \vspace{-0.5cm}
    \caption{Runtime comparison of all power flow models across all piecewise linear formulations.  The horizontal lines indicate the mean runtime of each piecewise linear formulation.}
    \label{fig:runtime_opf}
\end{figure}

\begin{table*}[th!]
\centering
\footnotesize
\caption{Quality and Runtime Results of AC Power Flow Formulations}
\begin{tabular}{|r|r|r||r|r|r|r||r|r|r|r|r|r|r|r|r|r|r|r|r|r|}
\hline
& & & \$/h & Delta & Delta & Delta & \multicolumn{4}{c|}{Runtime (seconds)} \\
Test Case & $|N|$ & $|E|$ & AC-$\lambda$ & AC-$\Delta$ & AC-$\Phi$ & AC-$\Psi$ & AC-$\lambda$ & AC-$\Delta$ & AC-$\Phi$ & AC-$\Psi$ \\
\hline
\hline
\multicolumn{11}{|c|}{Typical Operating Conditions (TYP)} \\
\hline
case\_500\_scenario\_140 & 500 & 733 & 3.0838e+05 & 0.00 & 0.00 & -0.00 & $<$1 & $<$1 & $<$1 & 4 \\
\hline
case\_500\_scenario\_33 & 500 & 733 & 2.9861e+05 & 0.00 & 0.00 & -0.00 & $<$1 & $<$1 & $<$1 & 4 \\
\hline
case\_500\_scenario\_72 & 500 & 733 & 4.6448e+05 & 0.00 & 0.00 & -0.00 & $<$1 & $<$1 & $<$1 & 6 \\
\hline
case\_793\_scenario\_12 & 793 & 913 & 2.6375e+05 & 0.00 & 0.00 & -0.00 & 2 & $<$1 & 2 & 3 \\
\hline
case\_793\_scenario\_135 & 793 & 913 & 2.1495e+05 & 0.00 & 0.00 & -0.00 & 2 & $<$1 & 2 & 2 \\
\hline
case\_793\_scenario\_72 & 793 & 913 & 2.6043e+05 & 0.00 & 0.00 & -0.00 & 2 & 2 & 2 & 3 \\
\hline
case\_2000\_scenario\_17 & 2000 & 3639 & 1.4495e+06 & -0.00 & 0.01 & -0.01 & 8 & 8 & 8 & 54 \\
\hline
case\_2000\_scenario\_65 & 2000 & 3639 & 9.3500e+05 & -0.00 & 0.00 & -0.01 & 8 & 8 & 8 & 77 \\
\hline
case\_2000\_scenario\_95 & 2000 & 3639 & 9.8239e+05 & -0.00 & 0.00 & -0.01 & 8 & 8 & 8 & 72 \\
\hline
case\_2312\_scenario\_172 & 2312 & 3013 & 3.8090e+05 & 0.00 & 0.00 & -0.01 & 6 & 6 & 7 & 8 \\
\hline
case\_2312\_scenario\_335 & 2312 & 3013 & 4.4250e+05 & 0.00 & 0.00 & -0.01 & 7 & 7 & 7 & 15 \\
\hline
case\_2312\_scenario\_60 & 2312 & 3013 & 4.1625e+05 & 0.00 & 0.00 & -0.01 & 7 & 7 & 7 & 10 \\
\hline
case\_2742\_scenario\_18 & 2742 & 4672 & 3.2300e+05 & -0.00 & -0.00 & -0.01 & 32 & 36 & 115 & 51 \\
\hline
case\_2742\_scenario\_36 & 2742 & 4672 & 3.5713e+05 & -0.01 & -0.00 & -0.01 & 31 & 35 & 39 & 49 \\
\hline
case\_2742\_scenario\_7 & 2742 & 4673 & 2.7608e+05 & -0.00 & -0.00 & -0.01 & 22 & 34 & 39 & 53 \\
\hline
case\_3022\_scenario\_136 & 3022 & 4135 & 4.9042e+05 & -0.00 & 0.00 & -0.01 & 9 & 10 & 10 & 10 \\
\hline
case\_3022\_scenario\_261 & 3022 & 4135 & 6.0265e+05 & -0.00 & -0.00 & -0.01 & 10 & 10 & 10 & 23 \\
\hline
case\_3022\_scenario\_45 & 3022 & 4135 & 5.3844e+05 & -0.00 & -0.00 & -0.01 & 10 & 9 & 10 & 12 \\
\hline
case\_3970\_scenario\_34 & 3970 & 6641 & 9.6241e+05 & -0.00 & 0.01 & -0.01 & 16 & 16 & 17 & 104 \\
\hline
case\_3970\_scenario\_13 & 3970 & 6643 & 6.3319e+05 & -0.00 & 0.01 & -0.02 & 15 & 18 & 21 & 93 \\
\hline
case\_3970\_scenario\_4 & 3970 & 6643 & 6.6357e+05 & -0.00 & 0.01 & -0.02 & 16 & 19 & 22 & 95 \\
\hline
case\_4020\_scenario\_102 & 4020 & 6986 & 6.2021e+05 & -0.01 & 0.00 & -0.02 & 18 & 23 & 26 & 61 \\
\hline
case\_4020\_scenario\_19 & 4020 & 6988 & 8.4275e+05 & -0.01 & 0.00 & -0.02 & 20 & 21 & 26 & 78 \\
\hline
case\_4020\_scenario\_47 & 4020 & 6988 & 8.2313e+05 & -0.00 & 0.00 & -0.02 & 14 & 18 & 20 & 93 \\
\hline
case\_4601\_scenario\_10 & 4601 & 7199 & 8.2755e+05 & 0.00 & 0.01 & -0.01 & 21 & 18 & 20 & 134 \\
\hline
case\_4601\_scenario\_3 & 4601 & 7199 & 8.6977e+05 & -0.00 & 0.01 & -0.01 & 21 & 19 & 20 & 135 \\
\hline
case\_4601\_scenario\_5 & 4601 & 7199 & 8.3184e+05 & 0.00 & 0.01 & -0.01 & 18 & 18 & 18 & 117 \\
\hline
case\_4619\_scenario\_17 & 4619 & 8150 & 4.7745e+05 & -0.01 & -0.01 & -0.02 & 18 & 22 & 22 & 64 \\
\hline
case\_4619\_scenario\_37 & 4619 & 8150 & 4.8192e+05 & -0.02 & -0.01 & -0.02 & 17 & 21 & 22 & 69 \\
\hline
case\_4619\_scenario\_4 & 4619 & 8152 & 7.9637e+05 & -0.01 & -0.00 & -0.02 & 19 & 21 & 26 & 85 \\
\hline
case\_4836\_scenario\_29 & 4836 & 7762 & 9.7137e+05 & 0.00 & 0.01 & -0.02 & 26 & 20 & 23 & 90 \\
\hline
case\_4837\_scenario\_13 & 4837 & 7765 & 7.9909e+05 & -0.01 & 0.00 & -0.02 & 17 & 19 & 20 & 82 \\
\hline
case\_4837\_scenario\_6 & 4837 & 7765 & 8.7327e+05 & -0.01 & 0.00 & -0.02 & 15 & 18 & 19 & 83 \\
\hline
case\_4917\_scenario\_167 & 4917 & 6726 & 1.3147e+06 & -0.00 & -0.00 & -0.01 & 18 & 17 & 18 & 38 \\
\hline
case\_4917\_scenario\_40 & 4917 & 6726 & 1.3911e+06 & -0.00 & 0.00 & -0.01 & 18 & 17 & 18 & 52 \\
\hline
case\_4917\_scenario\_99 & 4917 & 6726 & 1.2480e+06 & 0.00 & 0.00 & -0.01 & 18 & 17 & 19 & 25 \\
\hline
case\_8718\_scenario\_17 & 8718 & 14825 & 9.1927e+05 & 0.00 & 0.01 & -0.02 & 56 & 50 & 61 & 232 \\
\hline
case\_8718\_scenario\_5 & 8718 & 14825 & 1.0510e+06 & -0.00 & 0.01 & -0.02 & 65 & 55 & 61 & 235 \\
\hline
case\_9591\_scenario\_33 & 9591 & 15915 & 1.0635e+06 & 0.00 & 0.01 & -0.02 & 63 & 49 & 65 & 283 \\
\hline
case\_10000\_scenario\_15 & 10000 & 13193 & 1.0562e+06 & -0.00 & 0.00 & -0.06 & 139 & 36 & 34 & 207 \\
\hline
case\_10000\_scenario\_46 & 10000 & 13193 & 7.5720e+05 & -0.00 & -0.02 & -0.05 & 172 & 114 & 40 & 177 \\
\hline
case\_10000\_scenario\_81 & 10000 & 13193 & 1.4131e+06 & -0.00 & 0.00 & -0.06 & 43 & 45 & 46 & 87 \\
\hline
case\_10480\_scenario\_27 & 10480 & 18556 & 2.1314e+06 & -0.01 & 0.01 & -0.05 & 56 & 69 & 69 & 336 \\
\hline
case\_10480\_scenario\_43 & 10480 & 18557 & 2.0793e+06 & -0.02 & 0.00 & -0.05 & 52 & 61 & 61 & 290 \\
\hline
case\_10480\_scenario\_9 & 10480 & 18559 & 2.3159e+06 & -0.02 & -0.00 & -0.05 & 61 & 73 & 79 & 328 \\
\hline
case\_18889\_scenario\_89 & 18889 & 33435 & 2.2953e+06 & -0.01 & -0.02 & -0.06 & 149 & 143 & 186 & 593 \\
\hline
case\_19402\_scenario\_25 & 19402 & 34702 & 2.0834e+06 & -0.04 & -0.03 & -0.06 & 165 & 163 & 200 & 621 \\
\hline
case\_19402\_scenario\_42 & 19402 & 34704 & 1.9772e+06 & -0.04 & -0.03 & -0.06 & 182 & 190 & 242 & 548 \\
\hline
case\_24464\_scenario\_54 & 24464 & 37808 & 3.2305e+06 & -0.04 & 0.03 & -0.07 & 137 & 146 & 184 & 556 \\
\hline
case\_24464\_scenario\_37 & 24464 & 37816 & 2.6343e+06 & -0.00 & 0.02 & -0.10 & 130 & 124 & 128 & 570 \\
\hline
case\_24465\_scenario\_9 & 24465 & 37809 & 3.3570e+06 & -0.02 & -0.01 & -0.10 & 127 & 120 & 134 & 615 \\
\hline
case\_30000\_scenario\_19 & 30000 & 35393 & 1.0386e+06 & 0.00 & -0.01 & -0.04 & 1045 & 247 & 238 & 433 \\
\hline
case\_30000\_scenario\_36 & 30000 & 35393 & 1.1937e+06 & -0.00 & -0.01 & 0.16 & 748 & 301 & 237 & 1418 \\
\hline
case\_30000\_scenario\_5 & 30000 & 35393 & 1.0386e+06 & 0.00 & -0.01 & -0.04 & 239 & 257 & 484 & 392 \\
\hline
\end{tabular}\\
\label{tbl:deltas_time_ac}
\end{table*}

\subsection{Linear and Nonlinear Runtime Trends}

Given the stability of solution quality across the piecewise linear formulations, runtime performance becomes the most important criteria for comparison.
In this analysis the metric of interest is the relative runtime increase of a given piecewise linear formulation over the best runtime across all of the formulations.  Specifically, we define the {\em runtime-ratio} as,
%
\begin{align}
\mbox{runtime-ratio}_m = \frac{\mbox{runtime}_m}{\displaystyle \min_{n \in \{ \Psi, \lambda, \Delta, \Phi \}} \mbox{runtime}_n}  \;\; \forall m \in \{ \Psi, \lambda, \Delta, \Phi \} \label{eq:rt_ratio}
\end{align}
%
Figure \ref{fig:runtime_opf} presents the runtime-ratio of each formulation broken down by power flow formulation.  In these figures the y-axis presents the runtime-ratio in a log scale and the x-axis orders the 54 OPF problems from smallest to largest.  

{\em DC-OPF Results:}
With a variety of outliers occurring in each formulation, there is no clear winner in this formulation.
However, on average the $\lambda$ and $\Delta$ formulations perform best, which is consistent with similar benchmarking studies from the UC literature \cite{knueven2020mixed}.
It is important to highlight that the range of the y-axis in this case goes up to 1.7, which indicates that a poor selection of piecewise linear formulation can result in a performance reduction of no more than 70\%.
These results provide further validation of the experiment design as they replicate well-known results from the literature. 

{\em AC-OPF Results:} 
These results provide a stark contrast to the DC-OPF study.
The first significant difference is that the $\Psi$ model is significantly worse than all other formulations.
It is 5 times slower on average and can be 10 times slower in the worst case; that is 500\% and 1000\% slower, respectively.
The second significant difference is that although the remaining three models have similar performance, they all have notable outliers that are more than 2 times slower than the best formulation considered.
On average, the $\Delta$ formulation appears to be the best, but still suffers from notable outliers. 
A case-by-case runtime breakdown for the AC-OPF model is presented in Table \ref{tbl:deltas_time_ac}, for further inspection of specific cases.
Overall, these results highlight a drastically increased sensitivity of the AC-OPF problem to piecewise linear formulations.

{\em SOC-OPF Results:} 
By and large these results are similar to the AC-OPF results, which suggests that the increased sensitivity of the AC-OPF problem to piecewise linear formulations is a feature of the interior point algorithm, rather than an issue of non-convexity.
Two notable differences in these results are the overall reduction of outliers and the $\lambda$ formulation becoming a clear winner in terms of performance.
Overall, these results present novel insights into the impact of piecewise linear formulations on nonlinear optimization and suggest that a detailed analysis of nonlinear optimization algorithm behavior is required to find a consistent-best formulation for piecewise linear generation costs.

\subsection{Convex Nonlinear Algorithm Comparison}

The results from the previous section highlight a stark distinction between the behavior of linear programming algorithms (i.e., DC-OPF) and general purpose nonlinear interior point algorithms (i.e., AC-OPF and SOC-OPF).
However, the SOC-OPF model is a convex quadratic nonlinear model, which can be solved by specialized nonlinear optimization algorithms, such as Quadratic Constrained Quadratic Programming (QCQP) solvers.
This experiment briefly explores the possible benefit of convex quadratic optimization algorithms by comparing the solution of the SOC-OPF model via general purpose NLP (i.e., Ipopt) and a more specialized QCQP solver (i.e., Gurobi).

The first observation is that there is a significant difference in each solver's reliability, which is highlighted by Table \ref{tbl:soc_feasiblity}.
In at least 8 of the 54 networks considered, the QCQP solver reports a numerical error.  Interestingly, there is significant variability in the QCQP solver's reliability in different piecewise linear formulations.  Putting these reliability issues aside, the remaining analysis focuses on the subset of OPF cases where both solvers report convergence to an optimal solution.  In those cases, both solvers find solutions of identical quality, up to the numerical tolerances of the optimization algorithms, which suggests a correct implementation of both approaches.

\begin{table}[!ht]
\caption{Solver Reliability on SOC-OPF Models}
\label{tbl:soc_feasiblity}
\centering
\begin{tabular}{|l|l||r|r|r|r|r||r|r|r|r|r|r|r|r|r|r|r|r|}
\hline
Algorithm & Solver & $\lambda$ & $\Delta$ & $\Phi$ & $\Psi$ \\
\hline
\hline
NLP & Ipopt \cite{Ipopt} & 54 & 54 & 54 & 54 \\
\hline
QCQP & Gurobi \cite{gurobi} & 41 & 46 & 46 & 28 \\
\hline
\end{tabular}
\end{table}

The second observation is a significant difference in performance of the two algorithms.  Figure \ref{fig:soc_opf} presents a side-by-side comparison of the nonlinear solver runtimes for the subset of cases that both can solve.  In the figure, points below the diagonal line indicate a performance increase for the NLP solver and points above the line indicate a performance increase for the QCQP solver.  These results indicate the QCQP algorithm brings increased performance to the $\Psi$ model, which has very poor performance in the NLP solver.  However, for the other piecewise linear formulations, the NLP solver has a consistent performance advantage.  Overall, these results suggest that general purpose nonlinear interior point algorithms remain the most reliable and performant solution to large-scale convex nonlinear OPF problems.

\begin{figure}[t]
    \begin{center}
    \includegraphics[width=0.48\textwidth]{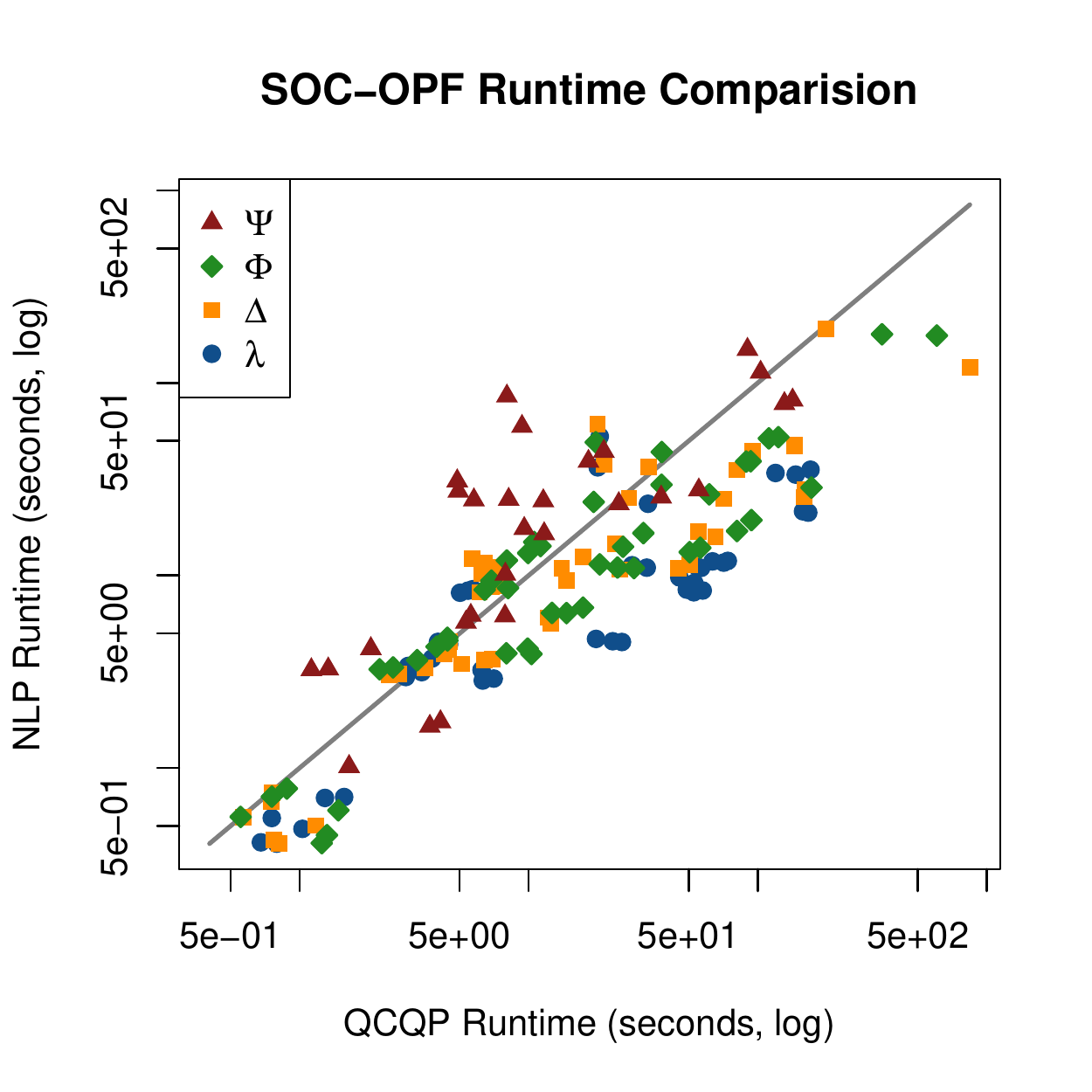}
    \end{center}
    \vspace{-0.5cm}
    \caption{A runtime comparison of solving the SOC-OPF problem with NLP and QCQP algorithms.}
    \label{fig:soc_opf}
\end{figure}

\subsection{Non-convex Nonlinear Solver Comparison}

The results from the previous sections highlight distinctions between the behavior of specialized convex optimization algorithms (e.g., linear programming and second order cone programming) and the general purpose nonlinear interior point algorithms, as implemented by Ipopt.
However, it is possible that the poor performance of the $\Psi$ model in the previous experiments is a consequence of a deficiency in Ipopt's implementation and is not inherent to interior point algorithms more broadly.
To better understand the consistency of the results presented in this work, the following experiment compares the performance of Ipopt to the commercial nonlinear optimization solver {\sc Knitro}.

Similar to the previous analysis, the first observation is that there is a difference in each solver's reliability, which is highlighted by Table \ref{tbl:ac_feasiblity}.
There appears to be a consistent convergence issue with two of the models considered, however this can likely be overcome by careful tuning of {\sc Knitro}'s convergence tolerance parameters.
Interestingly, the most significant variability in {\sc Knitro}'s reliability occurs in the $\Psi$ model, which provides further evidence that this formulation is particularly problematic for interior point algorithms.
Putting these reliability issues aside, the remaining analysis focuses on the subset of OPF cases where both solvers converge to a locally optimal solution.  In these cases, both solvers find solutions of identical quality, up to the numerical tolerances of the optimization algorithms, which suggests a correct implementation of both approaches.

\begin{table}[!ht]
\caption{Solver Reliability on AC-OPF Models}
\label{tbl:ac_feasiblity}
\centering
\begin{tabular}{|l|l||r|r|r|r|r||r|r|r|r|r|r|r|r|r|r|r|r|}
\hline
Algorithm & Solver & $\lambda$ & $\Delta$ & $\Phi$ & $\Psi$ \\
\hline
\hline
NLP & Ipopt \cite{Ipopt} & 54 & 54 & 54 & 54 \\
\hline
NLP & {\sc Knitro} \cite{Byrd2006} & 54 & 53 & 53 & 45 \\
\hline
\end{tabular}
\end{table}

The second observation are the trends in runtime performance of the two algorithms.
Figure \ref{fig:soc_opf} presents a side-by-side comparison of the solver runtimes for the subset of cases that both can solve.
In the figure, points below the diagonal line indicate a performance increase for Ipopt and points above the line indicate a performance increase for {\sc Knitro}.
These results indicate the two algorithms have similar performance on all of the models considered, which indicates that the performance challenges of the $\Psi$ model may persists in a variety of interior point algorithms.
It is worth noting that {\sc Knitro} tends to show slightly better performance on the $\Psi$ model while Ipopt has slightly better performance on the other models.
{\sc Knitro}'s implementation reduces the typical runtime of the $\Psi$ model by about half.  This is a notable improvement but it is not sufficient to change the overall conclusion of this work that the $\Delta$ and $\Lambda$ formulations a preferable for solving large-scale non-convex OPF problems.

\begin{figure}[t]
    \begin{center}
    \includegraphics[width=0.48\textwidth]{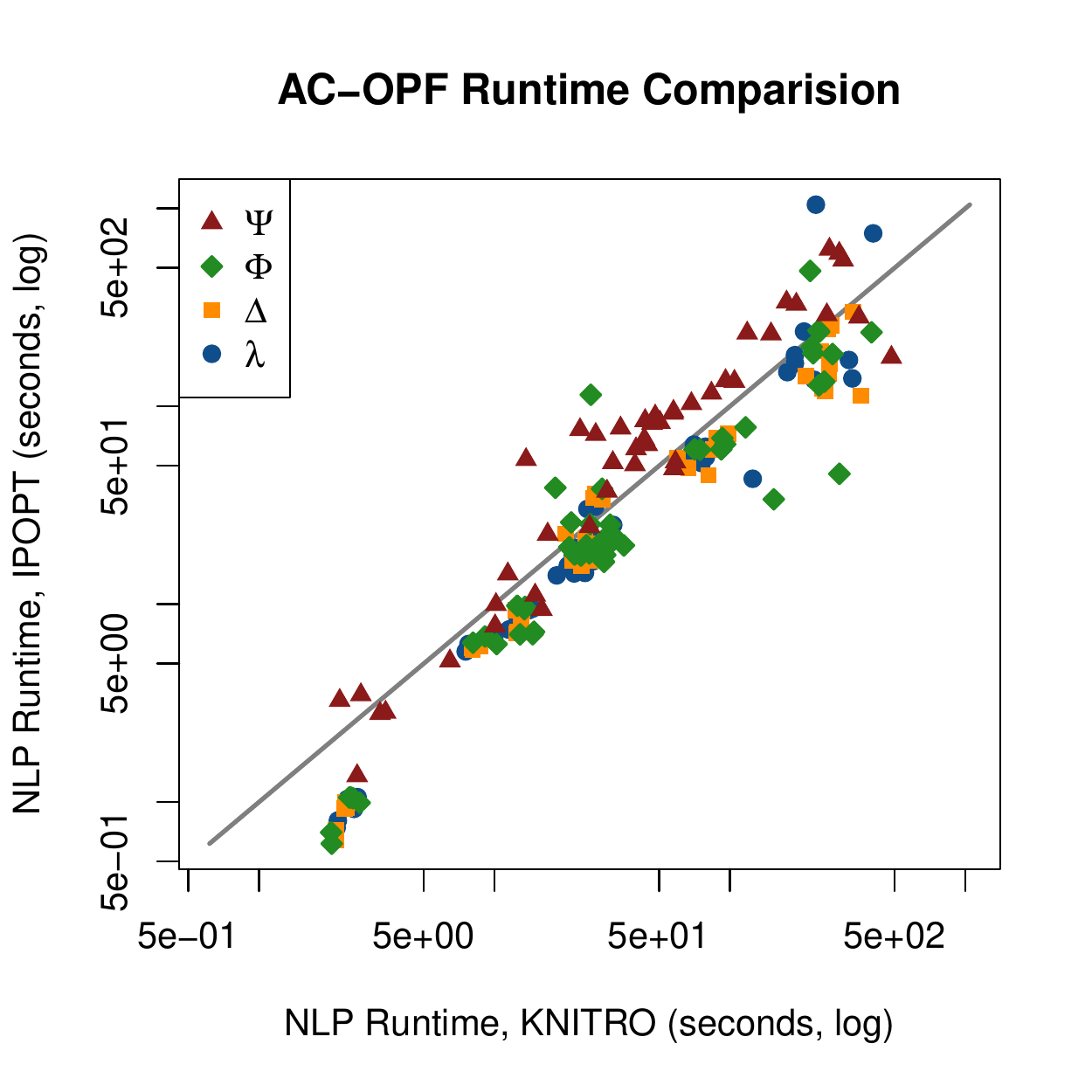}
    \end{center}
    \vspace{-0.5cm}
    \caption{A runtime comparison of solving the AC-OPF problem with Ipopt and {\sc Knitro} solvers.}
    \label{fig:soc_opf}
\end{figure}


\section{Conclusion}
\label{sec:conclusion}

Pursuing future nonlinear AC market clearing optimization algorithms, this work considered how to best formulate AC optimal power flow problems with piecewise linear generation cost functions.  To that end, core insights from the unit commitment literature in piecewise linear cost formulations \cite{knueven2020mixed} were combined with insights from the optimal power flow literature \cite{7271127,pm_pscc}, resulting in 12 variants of the optimal power flow problem.  A comprehensive numerical evaluation of these models on 54 realistic power network cases indicates that the ``$\lambda$'' and ``$\Delta$'' formulations of the piecewise cost functions prove to be the most suitable for current nonlinear optimization software, with the ``$\lambda$'' formulation being particularly suitable for convex relaxations of the power flow equations.  However, notable outliers in both models suggest ongoing research is required to ensure performance reliability of nonlinear optimization software, in preparation for real-world deployments with strict runtime requirements.

\section*{Acknowledgment}

The authors would like to thank Richard O'Neill for helpful feedback on a preliminary draft of this work.

\ifCLASSOPTIONcaptionsoff
  \newpage
\fi

\bibliographystyle{IEEEtran}
\bibliography{main}

LA-UR-20-23777






\end{document}